\documentclass[a4paper,11pt,twoside]{article}

\usepackage[T1]{fontenc}
\usepackage{amsfonts}
\usepackage{amssymb}
\usepackage{amsmath}

\usepackage{graphicx}
\usepackage{subfigure}

\usepackage{latexsym}

\usepackage{verbatim}

\makeatletter
\renewcommand{\@evenfoot}{\hfil - \thepage\ - \hfil}
\renewcommand{\@oddfoot}{\hfil - \thepage\ - \hfil}
\makeatother

\newcounter{theo}
\newcounter{propo}

\newenvironment{thm}
{\noindent\hangafter=1\hangindent=15pt\refstepcounter{theo}\textsc{Theorem \thetheo. --}\begin{sffamily}}
{\end{sffamily}\par}

\newenvironment{thm*}
{\noindent\hangafter=1\hangindent=15pt\refstepcounter{theo}\textsc{Theorem --}\begin{sffamily}}
{\end{sffamily}\par}

\newenvironment{thmt*}[1]
{\noindent\hangafter=1\hangindent=15pt\refstepcounter{theo}\textsc{Theorem -- #1}\\\begin{sffamily}}
{\end{sffamily}\par}

\newenvironment{prop}
{\noindent\hangafter=1\hangindent=15pt\refstepcounter{propo}\textsc{Proposition \thepropo. --}\begin{sffamily}}
{\end{sffamily}\par}

\newenvironment{prop*}
{\noindent\hangafter=1\hangindent=15pt\textsc{Proposition. --}\begin{sffamily}}
{\end{sffamily}\par}

\newenvironment{lemme}
{\noindent\hangafter=1\hangindent=15pt\refstepcounter{theo}\textsc{Lemma \thetheo. -- }\begin{sffamily}}
{\end{sffamily}\par}

\newenvironment{lemme*}
{\noindent\hangafter=1\hangindent=15pt\textsc{Lemma. -- }\begin{sffamily}}
{\end{sffamily}\par}

\newenvironment{coroll*}
{\noindent\hangafter=1\hangindent=15pt\textsc{Corollary. -- }\begin{sffamily}}
{\end{sffamily}\par}

\newenvironment{proof}
{\noindent\textsc{Proof.}\setlength{\parskip}{2pt}\setlength{\parindent}{0pt}}
{\hfill$\square$\par}

\newcommand{\bP}{\mathbb{P}}
\newcommand{\bE}{\mathbb{E}}

\newcommand{\eqd}{\begin{eqnarray*}} 
\newcommand{\eqf}{\end{eqnarray*}} 

\usepackage{hyperref}

\newcommand{\alphach}{{\check{\alpha}}}
\newcommand{\gammach}{{\check{\gamma}}}
\newcommand{\sigmach}{{\check{\sigma}}}
\newcommand{\ech}{{\check{e}}}
\newcommand{\Ech}{{\check{E}}}
\newcommand{\pch}{{\check{p}}}
\newcommand{\Gch}{{\check{G}}}

\begin{document}

\begin{center}
\Large Reversed Dirichlet environment and\\
Directional transience of random walks in Dirichlet random environment

\large
\textsc{Christophe Sabot and Laurent Tournier\footnote{
Universit\'e de Lyon ;
Universit\'e Lyon 1 ;
INSA de Lyon, F-69621 ;
Ecole Centrale de Lyon ;
CNRS, UMR5208, Institut Camille Jordan,
43 bld du 11 novembre 1918,
F-69622 Villeurbanne-Cedex, France\\
\emph{E-mail: }{\tt sabot@math.univ-lyon1.fr}, {\tt tournier@math.univ-lyon1.fr}\\This research was supported by the french ANR Project MEMEMO}}
\end{center}

\small
\paragraph{Abstract} We consider random walks in a random environment that is given by i.i.d.~Dirichlet distributions at each vertex of $\mathbb{Z}^d$ or, equivalently, oriented edge reinforced random walks on $\mathbb{Z}^d$. The parameters of the distribution are a $2d$-uplet of positive real numbers indexed by the unit vectors of $\mathbb{Z}^d$. We prove that, as soon as these weights are nonsymmetric, this random walk $(X_n)_n$ is transient in a direction (i.e.~it satisfies $X_n\cdot\ell\to_n+\infty$ for some $\ell$) with positive probability. In dimension 2, this result can be strenghened to an almost sure directional transience thanks to the 0-1 law from~\cite{Zerner-Merkl}. Our proof relies on the property of stability of Dirichlet environment by time reversal proved in~\cite{Sabot09}. In a first part of this paper, we also give a probabilistic proof of this property as an alternative to the change of variable computation used in~\cite{Sabot09}. 
\normalsize
\vspace{-.5cm}
\section*{Introduction}

\vspace{-.4cm}
Random walks in a multidimensional random environment have attracted considerable interest in the last few years. Unlike the one-dimensional setting, this model remains however rather little understood. Recent advances focused especially in two directions: small perturbations of a deterministic environment (cf.~for instance \cite{Bolthausen-Zeitouni}) and ballisticity (cf.~\cite{Sznitman} for a survey), but few conditions are completely explicit or known to be sharp. 


The interest for random Dirichlet environment stems from the desire to take advantage of the features of a \emph{specific} multidimensional environment distribution that make a few computations explicitly possible, and hopefully provide through them a better intuition for the general situation. The choice of Dirichlet distribution appears as a natural one considering its close relationship with linearly reinforced random walk on oriented edges (cf.~\cite{Pemantle} and~\cite{EnriquezSabot02}). The opportunity of this choice was further confirmed by the derivation of a ballisticity criterion by N.Enriquez and C.Sabot~\cite{EnriquezSabot06} (later improved by L.Tournier~\cite{Tournier09}), and more especially by the recent proof by C.Sabot~\cite{Sabot09} that random walks in Dirichlet environment on $\mathbb{Z}^d$ are transient when $d\geq 3$. 

The proof of the ballisticity criterion in~\cite{EnriquezSabot02} relies on an integration by part formula for Dirichlet distribution. This formula provides algebraic relations involving the Green function that allow to show that Kalikow's criterion applies. As for the proof of~\cite{Sabot09}, one of its key tools is the striking property (Lemma 1 of~\cite{Sabot09}) that, provided the edge weights have null divergence, a reversed Dirichlet environment (on a finite graph) still is a Dirichlet environment. 

In this paper we prove with this same latter property that random walks in Dirichlet environment on $\mathbb{Z}^d$ ($d\geq 1$) are directionally transient as soon as the weights are nonsymmetric. More precisely, under this condition, our result (Theorem~\ref{thm:transience}) states that directional transience happens with positive probability: for some direction $\ell$, $P_o(X_n\cdot\ell\to+\infty)>0$. Combined with S.Kalikow's 0-1 law from~\cite{Kalikow}, this proves that, almost surely, $|X_n\cdot \ell|\to+\infty$. Furthermore, in dimension 2, the 0-1 law proved by M.Zerner and F.Merkl (cf.~\cite{Zerner-Merkl}) enables to conclude that, almost surely, $X_n\cdot\ell\to+\infty$. 

The above mentioned property of the reversed Dirichlet environment was derived in~\cite{Sabot09} by means of a complicated change of variable. In Section~\ref{sec:reversal} of the present paper, we provide an alternative probabilistic proof of this important property. 

\section{Definitions and statement of the results}

Let us precise the setting in this paper. Let $G=(V,E)$ be a directed graph, i.e.~$V$ is the set of vertices, and $E\subset V\times V$ is the set of edges. If $e=(x,y)$ is an edge, we respectively denote by $\underline{e}=x$ and $\overline{e}=y$ its tail and head. We suppose that each vertex has finite degree, i.e.~that finitely many edges exit any vertex. $G$ is said to be strongly connected if for any pair of vertices $(x,y)$ there is a directed path from $x$ to $y$ in $G$. The set of environments on $G$ is the set
$$\Delta=\{(p_e)_{e\in E}\in(0,1)^E\ |\ \mbox{for all } x\in V, \sum_{e\in E,\underline{e}=x}p_e=1\}. $$
Let $(p_e)_{e\in E}$ be an environment. If $e=(x,y)$ is an edge, we shall write $p_e=p(x,y)$. Note that $p$ extends naturally to a measure on the set of paths: if $\gamma=(x_1,x_2,\ldots,x_n)\in V^n$ is a path in $G$, i.e.~if $(x_i,x_{i+1})\in E$ for $i=1,\ldots,n-1$, then we let $$p(\{\gamma\})=p(\gamma)=p(x_1,x_2)p(x_2,x_3)\cdots p(x_{n-1},x_n).$$
To any environment $p=(p_e)_{e\in E}$ and any vertex $x_0$ we associate the distribution $P^{(p)}_{x_0}$ of the Markov chain $(X_n)_{n\geq 0}$ on $V$ starting at $x_0$ with transition probabilities given by $p$. For every path $\gamma=(x_0,x_1,\ldots,x_n)$ in $G$ starting at $x_0$, we have $$P^{(p)}_{x_0}(X_0=x_0,\ldots,X_n=x_n)=p(\gamma).$$
When necessary, we shall specify by a superscript $P^{(p),G}_{x_0}$ which graph we consider.

Let $(\alpha_e)_{e\in E}$ be positive weights on the edges. For any vertex $x$, we let $$\alpha_x=\sum_{e\in E,\underline{e}=x} \alpha_e$$ be the sum of the weights of the edges exiting from $x$. The Dirichlet distribution with parameter $(\alpha_e)_{e\in E}$ is the law $\bP^{(\alpha)}$ of the random variable $p=(p_e)_{e\in E}$ in $\Delta$ such that the transition probabilities $(p_e)_{e\in E,\underline{e}=x}$ at each site $x$ are independent Dirichlet random variables with parameter $(\alpha_e)_{e\in E,\underline{e}=x}$. Thus, 
$$d\bP^{(\alpha)}(p)=\frac{\prod_{x\in V}\Gamma(\alpha_x)}{\prod_{e\in E}\Gamma(\alpha_e)}\prod_{e\in E}p_e^{\alpha_e-1}\prod_{e\in\widetilde{E}}dp_e,$$
where $\widetilde{E}$ is obtained from $E$ by removing arbitrarily, for each vertex $x$, one edge with origin $x$ (the distribution does not depend on this choice). The corresponding expectation will be denoted by $\bE^{(\alpha)}$. We may thus consider the probability measure $\bE^{(\alpha)}[P^{(p)}_{x_0}(\cdot)]$ on random walks on $G$. 

\paragraph{Time reversal} Let us suppose that $G$ is finite and strongly connected. Let $\Gch=(V,\Ech)$ be the graph obtained from $G$ by reversing all the edges, i.e.~by swapping the head and tail of the edges. To any path $\gamma=(x_1,x_2,\ldots,x_n)$ in $G$ we can associate the reversed path $\gammach=(x_n,\ldots,x_2,x_1)$ in $\Gch$. 

For any environment $p=(p_e)_{e\in E}$, the Markov chain with transition probabilities given by $p$ is irreducible, hence ($G$ being finite) we may define its unique invariant probability $(\pi_x)_{x\in V}$ on $V$ and the environment $\pch=(\pch_e)_{e\in\Ech}$ obtained by time reversal: for every $e\in E$,
$$\pch_\ech=\frac{\pi_{\underline{e}}}{\pi_{\overline{e}}}p_e.$$ 

For any family of weights $\alpha=(\alpha_e)_{e\in E}$, we define the reversed weights $\alphach=(\alphach_e)_{e\in\Ech}$ simply by $\alphach_\ech=\alpha_e$ for any edge $e\in E$, and the divergence of $\alpha$ is given, for $x\in V$, by
${\rm div}(\alpha)(x)=\alpha_x-\alphach_x.$

The following proposition is Lemma 1 from~\cite{Sabot09}. 

\begin{prop} \label{prop:reversal}
Suppose that the weights $\alpha=(\alpha_e)_{e\in E}$ have null divergence, i.e. $$\forall x\in V, {\rm div}(\alpha)(x)=0 \mbox{ (or equivalently }\alphach_x=\alpha_x).$$ Then, under $\bP^{(\alpha)}$, $\pch$ is a Dirichlet environment on $\Gch$ with parameters $\alphach=(\alphach_e)_{e\in \Ech}$. 
\end{prop}

In section~\ref{sec:reversal}, we give a neat probabilistic proof of this proposition. It also lies at the core of the proof of the directional transience. 

\paragraph{Directional transience on $\mathbb{Z}^d$} We consider the case $V=\mathbb{Z}^d$ ($d\geq 1$) with edges to the nearest neighbours. An i.i.d.~Dirichlet distribution on $G$ is determined by a $2d$-vector $\vec{\alpha}=(\alpha_1,\beta_1,\ldots,\alpha_d,\beta_d)$ of positive weights. We define indeed the parameters $\alpha$ on $E$ by $\alpha_{(x,x+e_i)}=\alpha_i$ and $\alpha_{(x,x-e_i)}=\beta_i$ for any vertex $x$ and index $i\in\{1,\ldots,d\}$. In this context, we let $P_o=\bE^{(\alpha)}[P^{(p)}_o(\cdot)]$ be the so-called annealed law. 

\begin{thm}\label{thm:transience}
Assume $\alpha_1>\beta_1$. Then $$\displaystyle P_o(X_n\cdot e_1\to_n +\infty)>0.$$
\end{thm}

The 0-1 law proved by S.Kalikow in Lemma 1.1 of~\cite{Kalikow} and generalized to the non-uniformly elliptic case by M.Zerner and F.Merkl in Proposition 3 of~\cite{Zerner-Merkl}, along with the 0-1 law proved by M.Zerner and F.Merkl in~\cite{Zerner-Merkl} (cf.~also~\cite{Zerner}) for the two-dimensional random walk allow to turn this theorem into almost-sure results: 

\begin{coroll*}\label{cor:corollaire}
Assume $\alpha_1>\beta_1$. 
\begin{enumerate}
	\item If $d\leq 2$, then $$P_o\mbox{-almost surely, }\hspace{.5cm}X_n\cdot e_1\to_n+\infty.$$
	\item If $d\geq 3$, then $$P_o\mbox{-almost surely, }\hspace{.5cm}|X_n\cdot e_1|\to_n+\infty.$$
\end{enumerate}
\end{coroll*}

\paragraph{Remark} This theorem provides examples of non-ballistic random walks that are transient in a direction. Proposition 12 of~\cite{Tournier09} shows indeed that the condition $2\sum_j(\alpha_j+\beta_j)-\alpha_i-\beta_i\leq1$ for a given $i$ implies $P_o$-almost surely $\tfrac{X_n}{n}\to_n 0$, because of the existence of ``finite traps'', so that any choice of small unbalanced weights (less than $\frac{1}{4d}$ for instance) is such an example. Another simple example where such a behaviour happens was communicated to us by A.Fribergh as well. It is however believed (cf.~\cite{Sznitman} p.227
) that transience in a direction implies ballisticity if a uniform ellipticity assumption is made, i.e.~if there exists $\kappa>0$ such that, for every edge $e$, $\bP$-almost surely $p_e\geq\kappa$. 

Under the hypothesis of Theorem~\ref{thm:transience} we conjecture that, for any $d\geq 1$, $$P_o(X_n\cdot e_1\to_n +\infty)=1,$$ and that the following identity is true: $$P_o(D=+\infty)=1-\frac{\beta_1}{\alpha_1},$$ where $D=\inf\{n\geq 0|X_n\cdot e_1<0\}$. In this paper, only the ``$\geq$'' inequality is proved. 

\section{Reversed Dirichlet environment. Proof of Proposition~\ref{prop:reversal}}\label{sec:reversal}

It suffices to prove that, for a given starting vertex $o$, the annealed distributions $\bE^{(\alpha)}[P_o^{(\pch)}(\cdot)]$ and $\bE^{(\alphach)}[P_o^{(p)}(\cdot)]$ on walks on $\Gch$ are the same. 

Indeed, the annealed distribution $\bE[P_o^{(p)}(\cdot)]$ characterizes the distribution $\bP$ of the environment, as soon as we assume that $\bP$ is a distribution on $\Omega=\prod_{x\in V}\mathcal{P}_x$ such that $\bP$-almost every environment $p$ is transitive recurrent. This results from considering the sample environment $p^{(n)}$ at time $n$:  
$$\mbox{for every $e=(x,y)\in E$, }\ \ p^{(n)}_e=\frac{\left|\{0\leq i < n| (X_i,X_{i+1})=e\}\right|}{\left|\{0\leq i <n| X_i=x\}\right|}.$$
From the recurrence assumption and the law of large numbers we deduce that $p^{(n)}$ makes sense for large $n$ and that, as $n\to\infty$, it converges almost-surely under $P_o^{(p)}$, and thus under $\bE[P_o^{(p)}(\cdot)]$, to a random variable $\widetilde{p}$ with distribution $\bP$. This way we can recover the distribution of the environment from that of the annealed random walk. 



We are thus reduced to proving that, for any finite path $\gamma$ in $G$, 
\begin{equation}\bE^{(\alpha)}[\pch(\gammach)]=\bE^{(\alphach)}[p(\gammach)].\label{eq:reverse}\end{equation} 

The only specific property of Dirichlet distribution to be used in the proof is the following ``cycle reversal'' property: 

\begin{lemme}\label{lem:reverse_cycle}
For any cycle $\sigma=(x_1,x_2,\ldots,x_n(=x_1))$ in $G$, $\bE^{(\alpha)}[p(\sigma)]=\bE^{(\alphach)}[p(\sigmach)]$, i.e. $$\bE^{(\alpha)}[p(x_1,x_2)\cdots p(x_{n-1},x_n)]=\bE^{(\alphach)}[p(x_n,x_{n-1})\cdots p(x_2,x_1)].$$ 
\end{lemme}

\begin{proof}
Remembering that the annealed random walk in Dirichlet environment is an oriented edge linearly reinforced random walk where the initial weights on the edges are the parameters of the Dirichlet distribution (cf.~\cite{EnriquezSabot02}), the left-hand side of the previous equality writes: 
$$\bE^{(\alpha)}[p(\sigma)]=\frac{\prod_{e\in E} \alpha_e(\alpha_e+1)\cdots (\alpha_e+n_e(\sigma)-1)}{\prod_{x\in V} \alpha_x(\alpha_x+1)\cdots (\alpha_x+n_x(\sigma)-1)},$$
where $n_e(\sigma)$ (resp.~$n_x(\sigma)$) is the number of crossings of the oriented edge $e$ (resp. the number of visits of the vertex $x$) in the path $\sigma$. The cyclicity gives $n_e(\sigma)=n_\ech(\sigmach)$ and $n_x(\sigma)=n_x(\sigmach)$ for all $e\in E,x\in V$. Furthermore, by assumption $\alphach_x=\alpha_x$ for every vertex $x$, and by definition $\alpha_e=\alphach_\ech$ for every edge $e$. This shows that the previous product matches the similar product with $\alphach$ and $\sigmach$ instead of $\alpha$ and $\sigma$, hence the lemma. 
\end{proof}

If $\sigma=(x_1,x_2,\ldots,x_n)$ is a cycle in $G$ (thus $x_n=x_1$), the definition of the reversed environment gives: 
\begin{eqnarray}
p(\sigma)
	& = & p(x_1,x_2)p(x_2,x_3)\cdots p(x_{n-1},x_n)\notag\\
	& = & \frac{\pi_{x_2}}{\pi_{x_1}}\pch(x_2,x_1)\frac{\pi_{x_3}}{\pi_{x_2}}\pch(x_3,x_2)\cdots \frac{\pi_{x_n}}{\pi_{x_{n-1}}} \pch(x_n,x_{n-1})\notag\\
	& = & \pch(x_2,x_1)\pch(x_3,x_2)\cdots\pch(x_n,x_{n-1})\notag\\
	& = & \pch(\sigmach)\label{eqn:reverse_cycle}
\end{eqnarray}
Therefore, the previous lemma gives: $\bE^{(\alpha)}[\pch(\sigmach)]=\bE^{(\alphach)}[p(\sigmach)]$, which is Equation~(\ref{eq:reverse}) for cycles. 

Consider now a non-cycling path $\gamma$ in $G$: $\gamma=(x=x_1,\ldots, x_n=y)$ with $x\neq y$. The same computation as above shows that $$\pch(\gammach)=\frac{\pi_x}{\pi_y}p(\gamma).$$
It is a well-known property of Markov chains that $\frac{\pi_x}{\pi_y}=E^{(p)}_y[N_x^y]$ where $N_x^y$ is the number of visits to $x$ before the next visit of $y$: $N_x^y=\sum_{i=0}^{H_y}{\bf 1}_{\{X_i=x\}}$, with $H_y=\inf\{n\geq 1\mid X_n=y\}$. It can therefore be written, in a hopefully self-explanatory schematical way: 
$$\frac{\pi_x}{\pi_y}=\sum_{k=0}^\infty P^{(p)}_y\left(N_x^y> k\right) = \sum_{k=0}^\infty p\left(\includegraphics[trim=0 26 0 26]{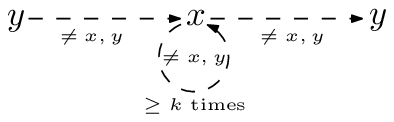}\right)=\sum_{k=0}^\infty p\left(\includegraphics[trim=0 25 0 25]{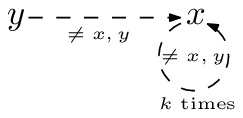}\right),$$

where the subscripts ``$\neq x,y$'' mean that the paths sketched by the dashed arrows avoid $x$ and $y$, and ``$(\geq) k$ times'' refers to the number of loops. Using the notation $\{\includegraphics[trim=0 2 0 2]{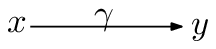}\}$ for the event where the walk follows the path $\gamma$ (which goes from $x$ to $y$), the Markov property at the time of $k$-th visit of $x$ gives: 
$$\pch(\gammach)=\sum_{k=0}^\infty p\left(\includegraphics[trim=0 25 0 25]{fig3c.eps}\right)p(\gamma)=\sum_{k=0}^\infty p\left(\includegraphics[trim=0 25 0 25]{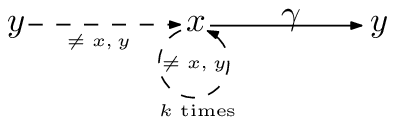}\right).$$

\vspace{.2cm}
The paths in the probability on the right-hand side are \emph{cycles}. Taking the expectation under $\bP^{(\alpha)}$ of both sides, we can use the Lemma to reverse them. We get: 
$$\bE^{(\alpha)}[\pch(\gammach)]=\sum_{k=0}^\infty  \bE^{(\alphach)}[p\left(\includegraphics[trim=0 25 0 25]{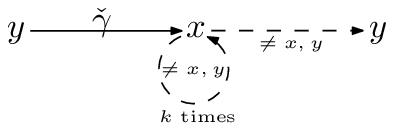}\right)].$$

It only remains to notice that the summation over $k\in\mathbb{N}$ allows to drop the condition on the path after it has followed $\gammach$:
$$\bE^{(\alpha)}[\pch(\gammach)]=\bE^{(\alphach)}[p(\gammach)].$$ 
This concludes the proof of the proposition. 

\section{Directional transience. Proof of Theorem~\ref{thm:transience}}

For any integer $M$, let us define the stopping times: 
$$T_M=\inf\{n\geq 0|X_n\cdot e_1\geq M\}$$
$$\widetilde{T}_M=\inf\{n\geq 0|X_n\cdot e_1\leq M\},$$
and in particular $D=\widetilde{T}_{-1}$. In addition, for any vertex $x$, we define $H_x=\inf\{n\geq 1\mid X_n=x\}$. 

Let $N,L\in\mathbb{N}^*$. Consider the finite and infinite ``cylinders'' $$C_{N,L}=\{0,\ldots,L-1\}\times(\mathbb{Z}_N)^{d-1} \subset C_N=\mathbb{Z}\times(\mathbb{Z}_N)^{d-1},$$ endowed with an i.i.d.~Dirichlet environment of parameter $\vec{\alpha}$. We are interested in the probability that the random walk starting at $o=(0,0,\ldots,0)$ exits $C_{N,L}$ to the right. Specifically, we shall prove: 
\begin{equation}\label{eqn:exit_cylinder}
	\bE^{(\alpha)}\left[P^{(p),C_N}_o(T_L<\widetilde{T}_{-1})\right]\geq 1-\frac{\beta_1}{\alpha_1}. 
\end{equation}


\begin{figure}[tbh]
	\centering
	\includegraphics[scale=2]{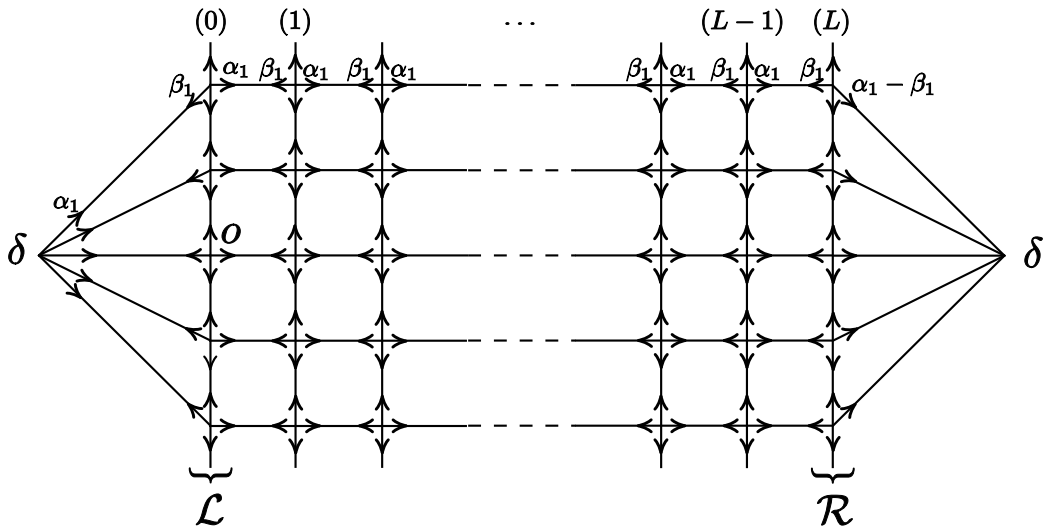}
	\caption{Definition of the weighted graph $G$. The top and bottom rows are identified (the graph is based on a cylinder) and the weights are the same on each line. Both vertices called $\delta$ are the same.}
	\label{fig:graph}
\end{figure}

To that aim, we introduce the finite weighted graph $G=G_{L,N}$ defined as follows (cf.~also figure~\ref{fig:graph}). The set of vertices of $G$ is $V=C_{N,L+1}\cup\{\delta\}$. Let us denote by $\mathcal{L}=\{0\}\times(\mathbb{Z}_N)^{d-1}$ and $\mathcal{R}=\{L\}\times(\mathbb{Z}_N)^{d-1}$ the left and right ends of the cylinder. The edges of $G$ inside $C_{N,L+1}$ go between nearest vertices for usual distance. In addition, ``exiting'' edges are introduced from the vertices in $\mathcal{L}$ and $\mathcal{R}$ into $\delta$, and ``entering'' edges are introduced from $\delta$ into the vertices in $\mathcal{L}$. The edges in $C_{N,L+1}$ (including exiting edges) are endowed with weights naturally given by $\vec{\alpha}$ \emph{except} that the edges exiting from $\mathcal{R}$ toward $\delta$ have weight $\alpha_1-\beta_1(>0)$ instead of $\alpha_1$. At last, the weight on the entering edges is $\alpha_1$. 

These weights, let us denote them by $\alpha'$, are chosen in such a way that ${\rm div }\, \alpha'(x)=0$ for every vertex $x$ (including $\delta$, necessarily). We are thus in position to apply Proposition~\ref{prop:reversal}.

Using Equation~(\ref{eqn:reverse_cycle}) to reverse all the cycles in $G$ that start at $\delta$ and get back to $\delta$ from $\mathcal{R}$ (for the first time), we get the following equality: 
$$P^{(p),G}_\delta\left(X_{H_\delta-1}\in \mathcal{R}\right) = P^{(\pch),\Gch}_\delta\left(X_1\in\mathcal{R}\right).$$
In the reversed graph $\Gch$, there are $N$ edges that exit $\delta$ to $\mathcal{L}$ (with reversed weight $\beta_1$), and $N$ edges that exit $\delta$ to $\mathcal{R}$ (with reversed weight $\alpha_1-\beta_1$). Combined with the Dirichlet distribution of $\pch$ under $\bP$ given by Proposition~\ref{prop:reversal} (notice that in fact we only reverse cycles, hence Lemma~\ref{lem:reverse_cycle} suffices), the previous equality gives: 
\begin{equation}\label{eqn:exit_G}
	\bE^{(\alpha')}\left[P^{(p),G}_\delta\left(X_{H_\delta-1}\in \mathcal{R}\right) \right] = \bE^{(\alphach')}\left[\sum_{e\in\Ech,\underline{e}=\delta}p(e)\right] =  \frac{N(\alpha_1-\beta_1)}{N(\alpha_1-\beta_1)+N\alpha_1}=1-\frac{\beta_1}{\alpha_1}.
\end{equation}
In the second equality we used the fact that under $\bP^{(\alphach')}$ the marginal variable $p_e$ follows a beta distribution with parameters $(\alphach'_e,(\sum_{e'\in\Ech,\underline{e'}=\delta}\alphach'_{e'})-\alphach'_e)$, and that the expected value of a Beta random variable with parameters $\alpha$ and $\beta$ is $\tfrac{\alpha}{\alpha+\beta}$. Alternatively, this equality follows also from the distribution of the first step of an oriented edge reinforced random walk with initial weights given by $\alpha'$. 

Let us show how this equality implies the lower bound~(\ref{eqn:exit_cylinder}). First, because of the symmetry of the weighted graph, the distribution of $X_1$ under $\bE^{(\alpha')}[P^{(p),G}_\delta(\cdot)]$ is uniform on $\mathcal{L}$, hence, since $o\in\mathcal{L}$: 
\begin{equation}\label{eqn:entering_point}
	\bE^{(\alpha')}\left[P^{(p),G}_\delta\left(X_{H_\delta-1}\in \mathcal{R}\right)\right]=\bE^{(\alpha')}\left[P^{(p),G}_o\left(X_{H_\delta-1}\in\mathcal{R}\right)\right].
\end{equation}
Then we relate the exit probability out of $C_{N,L}$ endowed with parameters $\alpha$ with the exit probability out of $C_{N,L+1}$ endowed with modified parameters $\alpha'$. These weights $\alpha$ and $\alpha'$ coincide in $C_{N,L}$; furthermore, a random walk on the cylinder $C_N$, from $o$, reaches the abcissa $L$ before the abcissa $L+1$, so that:
\begin{equation} \label{eqn:delete_column}
	\bE^{(\alpha)}\left[P^{(p),C_N}_o(T_L<\widetilde{T}_{-1})\right]\geq \bE^{(\alpha')}\left[P^{(p),C_N}_o(T_{L+1}<\widetilde{T}_{-1})\right]=\bE^{(\alpha')}\left[P^{(p),G}_o\left(X_{H_\delta-1}\in\mathcal{R}\right)\right].
\end{equation}
The equalities~(\ref{eqn:exit_G}) and~(\ref{eqn:entering_point}), along with~(\ref{eqn:delete_column}), finally give~(\ref{eqn:exit_cylinder}). 

To conclude the proof of the theorem, let us first take the limit as $N\to\infty$. Due to the ellipticity of Dirichlet environments, Lemma 4 of~\cite{Zerner-Merkl} shows that the random walk can't stay forever inside an infinite slab:
$$0=\bE^{(\alpha)}[P^{(p)}_o(T_L=\widetilde{T}_{-1}=\infty)]=\lim_{N\to\infty} \bE^{(\alpha)}[P^{(p)}_o(T_L\wedge \widetilde{T}_{-1}>T^\perp_N)]$$
where $T^\perp_N=\inf\{n\geq 0|\,|X_n-(X_n\cdot e_1)e_1|>N\}$ is the first time when the random walk is at distance greater than $N$ from the axis $\mathbb{R}e_1$. Using this with the stopping time $T^\perp_{N/2}$, we can switch from the cylinder $C_{N,L}$ to an infinite slab: 
\begin{eqnarray*}
\bE^{(\alpha)}[P^{(p),C_N}_o(T_L<\widetilde{T}_{-1})]
	& = & \bE^{(\alpha)}[P^{(p),C_N}_o(T_L<\widetilde{T}_{-1},T^\perp_{N/2}>T_L)]+\bE^{(\alpha)}[P^{(p),C_N}_o(T_L<\widetilde{T}_{-1},T^\perp_{N/2}\leq T_L)]\\
	& = & \bE^{(\alpha)}[P^{(p)}_o(T_L<\widetilde{T}_{-1},T^\perp_{N/2}>T_L)] + o_N(1)\\
	& = & \bE^{(\alpha)}[P^{(p)}_o(T_L<\widetilde{T}_{-1})] + o_N(1)',
\end{eqnarray*}
hence, with~(\ref{eqn:exit_cylinder}), $$\bE^{(\alpha)}[P^{(p)}_o(T_L<\widetilde{T}_{-1})]\geq 1-\frac{\beta_1}{\alpha_1}.$$ 

	
The limit when $L\to\infty$ simply involves decreasing events: 
$$P_o(D=+\infty)=\bE^{(\alpha)}[P_o^{(p)}(D=+\infty)]=\lim_{L\to\infty}\bE^{(\alpha)}[P^{(p)}_o(T_L<\widetilde{T}_{-1})]\geq 1-\frac{\beta_1}{\alpha_1}>0.$$
Moreover, Lemma 4 from~\cite{Zerner-Merkl} shows more precisely that if a slab is visited infinitely often, then \emph{both} half-spaces next to it are visited infinitely often as well. The event $\{D=\infty\}\cap\{\liminf_n X_n\cdot e_1<M\}$ has therefore null $P_o$-probability for any $M>0$, hence finally: 
$$P_o(X_n\cdot e_1\to_n+\infty)\geq \lim_{M\to\infty} P_o(D=\infty, \liminf_n X_n\cdot e_1\geq M)=P_o(D=\infty)>0.$$
This concludes the proof of the theorem.

\end{document}